\documentclass[12pt,a4paper]{article}
\usepackage{authblk}
\usepackage{indentfirst}
\usepackage{graphicx}
\usepackage{epsfig}

\usepackage{latexsym}
\usepackage{amsmath}
\usepackage{amssymb}
\usepackage{amsfonts}
\usepackage{mathrsfs}
\usepackage{pifont}

\usepackage{amsthm}
\usepackage{float}
\usepackage{subfigure}
\usepackage{color}
\usepackage{multicol}

\usepackage[colorlinks,linkcolor=blue,filecolor=black,citecolor=blue]{hyperref}

\newtheorem{thm}{Theorem}[section]
\newtheorem{lem}[thm]{Lemma}

\newtheorem{defi}{Definition}[section]
\newtheorem{coro}{Corollary}[section]
\newtheorem{rem}{\indent \bf Remark}[section]

\title{Unique Continuation on Quadratic Curves for Harmonic Functions}

\author[1]{Yufei Ke}
\author[2]{Yu Chen\thanks{Corresponding author: yuchen@sufe.edu.cn}}
\affil[1]{{\small School of Mathematics, Shanghai University of Finance and Economics, Shanghai, 200433, China}}
\affil[2]{{\small School of Mathematics, Shanghai University of Finance and Economics,
Shanghai 200433, China}}

\date{}
\begin{document}
\maketitle

\section*{Abstract}

 The unique continuation on quadratic curves for harmonic functions is discussed in this paper. By using complex extension method, the conditional stability of unique continuation along quadratic curves for harmonic functions is illustrated. The numerical algorithm is provided based on collocation method and Tikhonov regularization. The stability estimates on parabolic and hyperbolic curves for harmonic functions are demonstrated by numerical examples respectively.


\section{Introduction}

In this paper, we consider a problem that how to get global information by using local given information of the solution, namely, a unique continuation problem. The unique continuation for an elliptic equation states that if there is a connected open set $\Omega\in\mathbb{R}^{n}$ in which the solution $u$ lies and $\Xi\subset\Omega$ is an open subset of $\Omega$, then if $u|_{\Xi}=0$, $u|_{\Omega}\equiv0$. In two dimensions, the unique continuation means an analytic function must be zero if its zero points have accumulated points. We will first consider a parabolic curve, which is an analytic submanifold of two dimensional space. The aim is to study unique continuation on a parabola for harmonic functions, that is, if the value of function $u$ is known on a small part of the parabolic curve, then how to get the value of $u$ on a larger piece of the parabola. This problem is non-trivial since no measurement on the boundary $\partial \Omega$ is given. Moreover, it is an ill-posed problem due to its instability, which causes difficulties in numerical computations. For other types of quadratic functions, similar results can be obtained by the methods in this paper.

In the understanding of the unique continuation, unique continuation properties of various PDE's have been proved with development of PDEs theory \cite{Wolff1993}. Recently, the unique continuation on analytic sub-manifold has also been studied. For instance, Cheng and Yamamoto \cite{bib:3} fixed the measure points on a line and proved the conditional stability on a line for two-dimensional harmonic functions by adopting the method of complex continuation. Based on their results, Lu et al. \cite{bib:12} derived the unique continuation on a line for Helmholtz equation and the numerical treatment was firstly considered in their paper. The conditional stability of line unique continuation can also be applied to the estimation of unknown boundary for harmonic functions \cite{Bukhgeim1999}. Besides, the unique continuation has also been applied to the studies of other inverse problems, like control theory and optimal design problems. For example, the stability of soft obstacles reconstruction in inverse scattering was obtained by Isakov \cite{bib:4} based on unique continuation. Cheng et al. \cite{bib:13} studied the unique continuation for elasticity systems, which can be applied in phase transformation and anomalous diffusion in heterogeneous medium.

The focus of this paper is to generalize the conclusions of Cheng and Yamamoto \cite{bib:3} on a line to quadratic curves. A conditional stability for the continuation problem will be obtained based on the complex extension method, the unique continuation property is then implied. The main idea is to construct a holomorphic function based on the integral of single-layer potential with Green's function.
The determination of the analytic domain of the complex extention is a key to this process, which differs from the case of unique continuation on a line. After obtaining conditional stability estimates, a deterministic regularization, Tikhonov regularization, can be used to deal with the ill-posedness and construct stable numerical method \cite{bib:14}. Considering the discontinuity of real measurements, refer to \cite{bib:8}, the collocation method can be used to make specific numerical applications.

\begin{figure}[!ht]
\centering
    \includegraphics[width=2.6in]{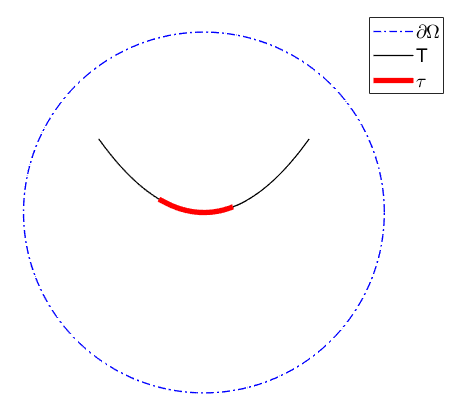}
\caption{Illustration of unique continuation on an quatratic curve $T$ with the measurement segment $\tau$ for a harmonic function $u$ in $\Omega$}
    \label{fig-illustration} 
\end{figure}

Our problem can be formulated as follows. Assume that $\Omega$ is a domain in $\mathbb{R}^{2}$, $\Gamma$, which is an analytic curve, is a continuous parabola in $\Omega$. Suppose that there exists open curves $\tau$ and $T$ satisfying: $\tau\subset T\subset\subset \Gamma$.
Consider a harmonic function $u(\mathbf{x})$ in $\Omega$, such that:
\begin{align}
  \Delta u(\mathbf{x})=0, \mathbf{x}\in\Omega.
\end{align}
Then this paper principally discusses the unique continuation on curve $T$ for $u$ which is known on $\tau$ (as illustrated in Fig. \ref{fig-illustration}), including the conditional stability estimate and the numerical computation. It should be pointed out that Cauchy values are not needed by our results, that is, no derivatives of harmonic functions on quadratic curves are required.

The rest of this paper is organized as follows: Section \ref{main} gives proof of the conditional stability on quadratic curves for harmonic functions, requiring that the taken curve segments have no intersections with boundaries. In Section \ref{example}, numerical applications implemented by the collocation method and Tikhonov's regularization are given. The numerical calculations on parabolic and logarithmic equations for harmonic functions, respectively, illustrate the validity of the conditional stability estimates. Conclusions are drawn in Section \ref{conclusion}.

\section{The unique continuation for harmonic functions}\label{main}

\begin{thm} Assume that $\tau$ and $T$ satisfy the previous definitions. Set $u(\mathbf{x})\in C^{2}(\Omega)$ as a harmonic function. If there exits a constant $M>0$, such that
\begin{align}\label{M}
  \| u\|_{C(\Omega)}\leq M,
\end{align}
then
\begin{align}\label{ux}
  \|u\|_{C(T)}\leq C\|u\|_{C(\tau)}^{\kappa},
\end{align}
where $C=C(M,\tau,T)$ is a positive constant, which is independent of $u$. $\kappa\in(0,1)$ depends on $\tau$ and $T$.
\end{thm}\label{thm}

We extend $u$ on a parabola from real plane to the complex plane first, which is a basis for the proof. Comparing with the extension of $u$ on a line, the determination of the analytic domain of the complex extention is a difficulty.

\begin{lem}\label{lemmab}\ \ Let $T_{1}\subset\mathbb{C}$ be the straight segment satisfing $T_{1}=\{(x,0)\mid (x,x^{2})\in T\}$, $B=\{z\in\mathbb{C}\mid dist(z,T_{1})<\varepsilon\}$, where $\varepsilon=\min\{\frac{\varepsilon_1}{4 diam(\Omega)},\frac{\varepsilon_1}{4}\}$, $\varepsilon_{1}=dist(T,\partial\Omega)$. Then $|(z-\zeta_{1})^{2}+(z^{2}-\zeta_{2})^{2}|\neq 0,\ \forall (\zeta_{1},\zeta_{2})\in\partial\Omega\subset\mathbb{R}^2,\ \forall z\in B\subset\mathbb{C}$.
\end{lem}
{\bf Proof}\ \ It can be proved by reduction to absurdity. It is obvious that $T_{1}\subset B$. For any $z\in B$, there exists a $z_{0}=x_{0}\in T_{1}$, such that $|z-z_{0}|<\varepsilon$.
If $\exists (\zeta_{1},\zeta_{2})\in\partial\Omega, z\in B$ such that $|(z-\zeta_{1})^{2}+(z^{2}-\zeta_{2})^{2}|= 0$, then
\begin{align}\label{z}
  \ \ &(z-\zeta_{1})^{2}+(z^{2}-\zeta_{2})^{2}=0 \notag\\
  \Rightarrow\ \ &z-\zeta_{1}=\pm\imath(z^{2}-\zeta_{2}) \notag\\
  \Rightarrow\ \ &(z-x_{0})+(x_{0}-\zeta_{1})=\pm\imath((z-x_{0}+x_{0})^{2}-\zeta_{2})=\pm\imath((z-x_{0})^{2}+(x_{0}^{2}-\zeta_{2})+2(z-x_{0})x_{0}) \notag\\
  \Rightarrow\ \ &|(x_{0}-\zeta_{1})\mp\imath(x_{0}^{2}-\zeta_{2})|=|\pm\imath((z-x_{0})^{2}+2(z-x_{0})x_{0})-(z-x_{0})|.
\end{align}
In (\ref{z}),we can see that
\begin{align*}
  LHS=|(x_{0}-\zeta_{1})\mp\imath(x_{0}^{2}-\zeta_{2})|\geq\varepsilon_{1},
\end{align*}
according to the definition of $\varepsilon_1$. On the other hand,
\begin{align*}
  RHS&=|\pm\imath((z-x_{0})^{2}+2(z-x_{0})x_{0})-(z-x_{0})| \\
     &\leq|z-x_{0}|^{2}+2|z-x_{0}||x_{0}|+|z-x_{0}| \\
     &\leq\varepsilon^{2}+2\varepsilon\cdot diam(\Omega)+\varepsilon \\
     &\leq\frac{\varepsilon_{1}}{4}\cdot\frac{\varepsilon_1}{4 diam(\Omega)}+\frac{3\varepsilon_1}{4} \\
     &<\varepsilon_{1},\ \
\end{align*}
where the last inequality is due to $\varepsilon_1\leq diam(\Omega)$. Then we arrives a contradiction and the proof is complete.

\begin{lem}\label{lemmaa}\ \ Let $\Omega\subset\mathbb{R}^{2}$, then for any harmonic function $u=u(x,y)$ in $\Omega$, which satisfies $\|u\|_{C(\Omega)}\leq M$, $M>0$ is a constant, there exist a simply connected domain $B\subset\mathbb{C}$ which contains $T_1$ defined in lemma 2.1, and a holomorphic function $v=v(z)$ ($z=\lambda + \imath \mu\in\mathbb{C}$) in $B$, such that
\begin{align}
  v(\lambda,\nu)=u(x,y),\ \ y=x^{2}, \lambda=x, \nu=0, (x,y)\in T.
\end{align}
\end{lem}

{\bf Proof}\ \ For any $\Omega^{\prime}$ satisfying $\Gamma\subset\Omega^{\prime}\subset\Omega$, since $u$ is harmonic in $C(\overline{\Omega^{\prime}})$, by Green's formula, there exists a density $\mu$ in $\partial \Omega^{\prime}$, such that
\begin{align}\label{uGreen}
  u(x,y)=\int_{\partial \Omega^{\prime}}\log((x-\zeta_{1})^{2}+(y-\zeta_{2})^{2})\mu(\mathbf{\zeta})ds_{\mathbf{\zeta}},\ \ \mathbf{x}=(x,y)\in \Omega^{\prime},
\end{align}
For the parabola $y=x^2$, we define the complex function,
\begin{align} \label{v}
  v(\lambda+\imath\nu)=\int_{\partial \Omega^{\prime}}\log((\lambda+\imath\nu-\zeta_{1})^{2}+((\lambda+\imath\nu)^{2}-\zeta_{2})^{2})\mu(\mathbf{\zeta})ds_{\mathbf{\zeta}},\quad \lambda+\imath \nu\in \mathbb{C}.
\end{align}
It can be seen that when $\lambda=x, \nu=0$,
\begin{align*}
  v(\lambda)=\int_{\partial \Omega^{\prime}}\log((\lambda-\zeta_{1})^{2}+(\lambda^{2}-\zeta_{2})^{2})\mu(\mathbf{\zeta})ds_{\mathbf{\zeta}}=u(\lambda,\lambda^2), \lambda\in\mathbb{R}.
\end{align*}

According to Lemma \ref{lemmab}, there exists a so called analytic domain $B\subset\mathbb{C}$ satisfying $T_{1}\subset B$, such that $\log((\lambda-\zeta_{1})^{2}+(\lambda^{2}-\zeta_{2})^{2})$ is analytic in $B$.

On the other hand, it was shown in \cite{bib:5} that when $u\in H^{2}(\partial \Omega^{\prime})$, $\mu\in H^{1}(\partial \Omega^{\prime})$, where $H^{s}(\partial \Omega^{\prime})(s=1,2)$ is Sobolev space on $\partial \Omega^{\prime}$. One has
\begin{align}
  \frac{1}{C_{1}}\|u\|_{H^{2}(\partial \Omega^{\prime})}\leq\|\mu\|_{H^{1}(\partial \Omega^{\prime})}\leq C_{1}\|u\|_{H^{2}(\partial \Omega^{\prime})},
\end{align}
here $C_{1}$ is a constant.
According to Sobolev embedding theorem and interior elliptic estimate \cite{Gilbarg1983}, there exists a constant $C_3=C_3(\Omega',\Omega)>0$, which is independent of $u$, such that
\begin{align}\label{G}
  \|\mu\|_{C(\partial \Omega^{\prime})}\leq C_2\|u\|_{C^2(\Omega')} \leq C_3\|u\|_{C(\Omega)}= C_3 M.
\end{align}
Then $v$ is analytic in $B$ and
\begin{align*}
  v(\lambda+\imath\nu)=u(x,y),\ \text{for} \ y=x^{2}, \lambda=x, \nu=0, (x,y)\in T.
\end{align*}

For other quadratic curves like a hyperbolic curve, if we take the observation range on one branch of the hyperbola, a similar method can be used to construct a complex extension of the harmonic function $u$. Without lose of generality, set $\Gamma=\{(x,y)\in\Omega\mid y^{2}-x^{2}=\frac{1}{9}, y>0\}$.

\begin{lem}\label{lemmad}\ \ Let $T_{1}\in\mathbb{C}$ be the straight segment satisfing $T_{1}=\{(x,0)\mid (x,\sqrt{x^{2}+\frac{1}{9}})\in T\}$, $B=\{z\in\mathbb{C}\mid dist(z,T_{1})<\min\{\frac{\varepsilon_1}{6},\frac{\varepsilon_1}{36 diam(\Omega)}\}\}$, where $\varepsilon_{1}=\min\{1,dist(T,\partial\Omega)\}$. Then $|(z-\zeta_{1})^{2}+(\sqrt{z^{2}+\frac{1}{9}}-\zeta_{2})^{2}|\neq 0, \forall (\zeta_{1},\zeta_{2})\in\partial\Omega,\ \forall z\in B\subset\mathbb{C}$.
\end{lem}

{\bf Proof}\ \ It can be proved by contradiction. It is obvious that $T_{1}\subset B$. For any $z\in B$, there exists a $z_{0}=x_{0}\in T_{1}$, such that $|z-z_{0}|<\varepsilon$.
If $\exists (\zeta_{1},\zeta_{2})\in\partial\Omega, z\in B$ such that $ |(z-\zeta_{1})^{2}+(\sqrt{z^{2}+\frac{1}{9}}-\zeta_{2})^{2}|= 0$, then
\begin{align}\label{z1}
  \ \ &(z-\zeta_{1})^{2}+(\sqrt{z^{2}+\frac{1}{9}}-\zeta_{2})^{2}=0 \notag\\
  \Rightarrow\ \ &z-\zeta_{1}=\pm\imath(\sqrt{z^{2}+\frac{1}{9}}-\zeta_{2}) \notag\\
  \Rightarrow\ \ &(z-x_{0})+(x_{0}-\zeta_{1})=\pm\imath(\sqrt{(z-x_{0}+x_{0})^{2}+\frac{1}{9}}-\zeta_{2}) \notag\\
  &\ \ \ \ \ =\pm\imath(\sqrt{(z-x_{0})^{2}+2(z-x_{0})x_{0}+(x_{0}^{2}+\frac{1}{9})}+(\sqrt{x_{0}^{2}+\frac{1}{9}}-\zeta_{2})-\sqrt{x_{0}^{2}+\frac{1}{9}}) \notag\\
  \Rightarrow\ \ &|(x_{0}-\zeta_{1})\mp\imath(\sqrt{x_{0}^{2}+\frac{1}{9}}-\zeta_{2})|\notag\\
&\ \ \ \ \ =|\pm\imath(\sqrt{(z-x_{0})^{2}+2(z-x_{0})x_{0}+(x_{0}^{2}+\frac{1}{9})}-\sqrt{x_{0}^{2}+\frac{1}{9}})-(z-x_{0})|.
\end{align}

From (\ref{z1}),we can see that
\begin{align*}
  LHS=|(x_{0}-\zeta_{1})\mp\imath(\sqrt{x_{0}^{2}+\frac{1}{9}}-\zeta_{2})|\geq dist(T,\partial \Omega)\geq\varepsilon_{1},
\end{align*}
and
\begin{align*}
  RHS&=|\pm\imath(\sqrt{(z-x_{0})^{2}+2(z-x_{0})x_{0}+(x_{0}^{2}+\frac{1}{9})}-\sqrt{x_{0}^{2}+\frac{1}{9}})-(z-x_{0})| \\
     &\leq|\sqrt{(z-x_{0})^{2}+2(z-x_{0})x_{0}+(x_{0}^{2}+\frac{1}{9})}-\sqrt{x_{0}^{2}+\frac{1}{9}}|+|z-x_{0}| \\
     &\leq \frac{|(z-x_0)^2+2(z-x_0)x_0|}{|\sqrt{(z-x_0)^2+2(z-x_0)x_0+(x_0^2+\frac{1}{9})}+\sqrt{x_0^2+\frac{1}{9}}|}+|z-x_0|.\\
     %
    %
\end{align*}
Since $|z-x_0|\leq \varepsilon$, with the definition of $\varepsilon$ and denote $|\Omega|=diam(\Omega)$ we have
\[
|(z-x_{0})^{2}+2(z-x_{0})x_{0}|\leq \varepsilon^2+2 |\Omega|\varepsilon\leq \frac{\varepsilon_1^2}{36}+2 |\Omega|\frac{\varepsilon_1}{36 |\Omega|}\leq \frac{\varepsilon_1\cdot 1}{36}+\frac{\varepsilon_1}{18}=\frac{\varepsilon_1}{12}.
\]
Due to $\frac{\varepsilon_1}{12}\leq \frac{1}{12}\leq \frac{1}{9}+x_0^2$, in the corresponding one-valued branch we have the real part
\[
\mathrm{Re} \left(\sqrt{(z-x_0)^2+2(z-x_0)x_0+(x_0^2+\frac{1}{9})}\right)\geq 0 .
\]
Therefore,
\[
RHS\leq \frac{\frac{\varepsilon_1}{12} }{\sqrt{x_0^2+\frac{1}{9}}}+\frac{\varepsilon_1}{6}< \frac{\varepsilon_1}{2}.
\]
Which is a contradiction with $LHS\geq \varepsilon_1$. Thus the proof is complete.

\begin{lem}\label{lemmac}\ \ Assume that $u=u(x,y)$ is harmonic in $\Omega\subset\mathbb{R}^{2}$, which satisfies $\|u\|_{C(\Omega)}\leq M$ where $M>0$ is a constant, then there exist a simply connected domain $B\subset \mathbb{C}$ containing the line $T_1$ defined in lemma 3, and a holomorphic function $v=v(\lambda,\nu)$ in $B$, such that
\begin{align}
  v(\lambda,\nu)=u(x,y),\ \ y^{2}-x^{2}=\frac{1}{9}, \lambda=x, \nu=0, (x,y)\in T.
\end{align}
\end{lem}

{\bf Proof}\ \ For any $\Omega^{\prime}$ which satisfies $\Gamma\subset\Omega^{\prime}\subset\Omega$, by Green's formula, since $u$ is harmonic in $C(\overline{\Omega^{\prime}})$, there exists a density $\mu$ in $\partial \Omega$, such that

\begin{align}\label{uGreenh}
  u(x,y)=\int_{\partial \Omega^{\prime}}\log((x-\zeta_{1})^{2}+(y-\zeta_{2})^{2})\mu(\mathbf{\zeta})ds_{\mathbf{\zeta}},\ \ \mathbf{x}=(x,y)\in \Omega^{\prime}.
\end{align}

Then according to Lemma \ref{lemmad} and similar to proof of Lemma \ref{lemmaa}, there exists a domain $B\in\mathbb{C}$, such that $v$ defined by
\begin{align} \label{vh}
  v(\lambda,\nu)=\int_{\partial \Omega^{\prime}}\log((\lambda+\imath\nu-\zeta_{1})^{2}+(\sqrt{(\lambda+\imath\nu)^{2}+\frac{1}{9}}-\zeta_{2})^{2})\mu(\mathbf{\zeta})ds_{\mathbf{\zeta}},\ \ z=\lambda+\imath\mu \in B,
\end{align}
is holomorphic in $B$. On $\mathbf{x}=(x,\sqrt{x^{2}+\frac{1}{9}})\in T\subset\Omega^{\prime}$,
\begin{align*}
  u(x,y)=\int_{\partial \Omega^{\prime}}\log((x-\zeta_{1})^{2}+(\sqrt{x^{2}+\frac{1}{9}}-\zeta_{2})^{2})\mu(\mathbf{\zeta})ds_{\mathbf{\zeta}},\ \ \mathbf{x}=(x,\sqrt{x^{2}+\frac{1}{9}})\in T.
\end{align*}
and on $\lambda=x, \nu=0$, one has
\begin{align*}
  v(\lambda,\nu)=\int_{\partial \Omega^{\prime}}\log((\lambda-\zeta_{1})^{2}+(\sqrt{\lambda^{2}+\frac{1}{9}}-\zeta_{2})^{2})\mu(\mathbf{\zeta})ds_{\mathbf{\zeta}}.
\end{align*}
which indicates,
\begin{align*}
  v(\lambda,\nu)=u(x,y),\ \ y^{2}-x^{2}=\frac{1}{9}, \lambda=x, \nu=0, (x,y)\in T.
\end{align*}

Denote $D$ the domain satisfying
\begin{align}
  D=\{\lambda+\imath \nu\in \mathbb{C} \mid (\lambda,\nu)\in [\lambda_{0},\lambda_{0}+r]\times[-h,h]\}.
\end{align}
Let $l=\{(x,0)\mid x\in[\lambda_0+r_{1},\lambda_0+r_{2}]\}\subset D$ be a closed segment, where $0<r_{1}<r_{2}<r$. A harmonic measure with respect to $D$ and $l$ can be defined as follows.

\begin{defi} \label{hm}{\rm (harmonic measure)}\ \ A function $\varphi(\lambda,\nu)$ is called a harmonic measure for $D$ and $l$, if $\varphi(\lambda,\nu)$ satisfies
\begin{align*}
  &\Delta\varphi(\lambda,\nu) = 0,\ \ (\lambda,\nu)\in D\setminus l, \\
  &\varphi(\lambda,\nu) = 0,\ \ (\lambda,\nu)\in\partial D,\\
  &\varphi(\lambda,\nu) = 1,\ \ (\lambda,\nu)\in l.
\end{align*}
\end{defi}

On the basis of Friedman and Vogelius\cite{bib:6} and Kellogg\cite{bib:7}, $\varphi(z)$ exists and is unique. $\varphi\in C^{\tau}(\overline{D})(0<\tau<1)$ can be obtained by \cite{bib:6}. For the estimate of the harmonic measure, we have

\begin{lem}\label{hmestimate}\ \ There exists a positive constant $C_{4}$, such that the harmonic measure defined in definition \ref{hm} satisfies
\begin{align}
  \varphi(x,0)\geq C_{4}(\lambda_0+r-x),\ \ x\in[\lambda_0+r_{2},\lambda_0+r].
\end{align}
\end{lem}
{\bf Proof}\ \ 
From definition \ref{hm} and the maximum principle for harmonic functions,
\begin{align}\label{ph}
  0<\varphi(\lambda,\nu)<1,\quad \lambda+\imath \nu\in D\backslash l. 
\end{align}
Assume to the contrary that $\forall \epsilon>0$ there exists $x_\epsilon\in  [\lambda_0+r_2,\lambda_0+r]$ such that $\varphi(x_\epsilon,0)<\epsilon (\lambda_0+r-x_\epsilon)$. Then, we can take $\epsilon=\frac{1}{n}, n=1,2,...$ and have the corresponding $x_n$ satisfying
\begin{align}\label{phixn}
 \varphi(x_{n},0) < \frac{1}{n}(\lambda_0+r-x_{n})\leq \frac{1}{n}(r-r_2)\longrightarrow 0,\ \ n\rightarrow\infty.
\end{align}
Note that if some $x_n=\lambda_0+r$ then it contradicts to the boundary condition since $(\lambda_0+r,0)\in \partial D$. It can also be seen that there are infinity many different $x_n$. Thus we can assume $x_n\neq\lambda_0+r$ and further have 
\begin{equation}\label{phixn2}
 \frac{\varphi(x_{n},0)}{\lambda_0+r-x_n} < \frac{1}{n} \longrightarrow 0,\ \ n\rightarrow\infty.
\end{equation}
Since $[\lambda_0+r_{2},\lambda_0+r]$ is compact, there must exist an $\hat{x}\in[\lambda_0+r_{2},\lambda_0+r]$, and a subsequence $\{x_{n_k}\}$ such that $x_{n_k}\longrightarrow\hat{x}$ for $k\rightarrow\infty$. If $\hat{x}\neq \lambda_0+r$ then $\varphi(\hat{x},0)=0$ due to \eqref{phixn} and the continuity of $\varphi$, which is a contradition to \eqref{ph}. If $\hat{x}= \lambda_0+r$ then from \eqref{phixn2} and the boundary geometry we have 
\begin{align*}
  \frac{\partial\varphi(\lambda,\nu)}{\partial\lambda}\Big|_{(\lambda_0+r,0)}=0,
\end{align*}
which is a contradiction to the maximum principle. Therefore, the conclusion of the lemma is true.

\begin{lem}\label{stability}\ \ Assume that $v=v(\lambda,\nu)$ is a holomorphic function in $D$. If $|v(\lambda,\nu)|\leq M_{1},(\lambda,\nu)\in D$, let $\varepsilon=\max_{(\lambda,\nu)\in l}|v(\lambda,\nu)|$,then
\begin{align}
  |v(x)|\leq M_{1}(\frac{\varepsilon}{M_{1}})^{\varphi(x)}\leq M_{1}(\frac{\varepsilon}{M_{1}})^{C_{5}(\lambda_0+r-x)},\ \ x\in[\lambda_0+r_{2},\lambda_0+r],
\end{align}
where $C_{5}$ is a constant.
\end{lem}
It can be proved by the same method in \cite{bib:3} with the conclusion of lemma \ref{hmestimate}.

So far, it is sufficient to start with the derivation of the main theorem.

{\bf Proof}\ \  Consider the parabolic curve first. Without lose of generality, suppose that
\begin{align}
  \Gamma = \{\mathbf{x}=(x,y)\mid y=x^{2},0<x<R\},
\end{align}
and
\begin{align*}
  &T = \{\mathbf{x}=(x,y)\mid y=x^{2},a<x<b\}, \\
  &\tau = \{\mathbf{x}=(x,y)\mid y=x^{2},c<x<d\},
\end{align*}
here $0<a<c<d<b<R$.
If $\Gamma$ is a  hyperbola, without lose of generality, suppose that
\begin{align}
  \Gamma = \{\mathbf{x}=(x,y)\mid y^{2}-x^{2}=\frac{1}{9},0<x<R\},
\end{align}
and
\begin{align*}
  &T = \{\mathbf{x}=(x,y)\mid y^{2}-x^{2}=\frac{1}{9},a<x<b\}, \\
  &\tau = \{\mathbf{x}=(x,y)\mid y^{2}-x^{2}=\frac{1}{9},c<x<d\},
\end{align*}
here $0<a<c<d<b<R$.

Accordingly, suppose that
\begin{align*}
  &\Gamma_{1} = \{\mathbf{x}=(x,0)\mid 0<x<R\}, \\
  &T_{1} = \{\mathbf{x}=(x,0)\mid a<x<b\}, \\
  &\tau_{1} = \{\mathbf{x}=(x,0)\mid c<x<d\}.
\end{align*}

Then by lemma \ref{lemmab} and \ref{lemmad} there exists the corresponding domain $B$ and subsequently a positive constant $\rho$ can be found such that $(a-\rho,b+\rho)\times(-\rho,\rho)\subset\overline{B}$. $v=v(\lambda,\nu)$ in lemma \ref{lemmaa} and lemma \ref{lemmac} is holomorphic in $(a-\rho,b+\rho)\times(-\rho,\rho)$. Due to the maximum principle,
\begin{align}
  \|v\|_{C((a-\rho,b+\rho)\times(-\rho,\rho))}\leq \|v\|_{C(B)}.
\end{align}
By (\ref{G}),
\begin{align}
  \|v\|_{C((a-\rho,b+\rho)\times(-\rho,\rho))}\leq \|v\|_{C(B)}\leq C_{6}\|u\|_{C(\Omega)}\leq C_{6}M.
\end{align}
where $C_{6}$ is a constant.

Then by Lemma \ref{stability}, there exist $C_{7}$ and $C_{8}$ just depending on $T$ and $\Omega$ such that,
\begin{align}
  |u(\mathbf{x})|=|v(x,0)|\leq C_{6}M(\frac{\varepsilon}{C_{6}M})^{C_{7}(b+\rho-x)}\leq C_{8}\varepsilon^{C_{7}\rho},\ \ x\in[d,b].
\end{align}
Similarly,
\begin{align}
  |u(\mathbf{x})|=|v(x,0)|\leq C_{6}M(\frac{\varepsilon}{C_{6}M})^{C_{9}(x+\rho-a)}\leq C_{10}\varepsilon^{C_{9}\rho},\ \ x\in[a,c].
\end{align}

\begin{coro}\ \ (\ref{ux}) is called a conditional stability estimate of $u$ with condition (\ref{M}), which implies the unique continuation, i.e. $u|_{\tau}=0$ yields $u|_{T}\equiv0$.
\end{coro}

\begin{rem}\ \ Replace the parabola with other quadratic curves, or even other analytical curves, the conditional stability estimate may be obtained similarly.
\end{rem}

\begin{rem}\ \ From the proof, it indicates that the degree of control for $u$ on $T$ is related to the harmonic measure $\varphi$, in other words, $\kappa$ in Theorem \ref{thm} may be further quantified by the harmonic measure $\varphi$ which is determined by $T$, $\tau$ and $\Omega$.
\end{rem}

\begin{rem}\ \ The unique continuation does not hold outside the quadratic curve, but holds on the quadratic curve only. However, if the curve is a high-order curve, like $y=x^{n}, n>2$, the unique continuation would hold outside the quadratic curve. For details, refer to \cite{bib:15}.
\end{rem}

\section{Numerical methods and applications}\label{example}

The numerical method to achieve unique continuation on a parabola for harmonic functions is given in this section.

Assume that $u$ is the solution of the harmonic equation below
\begin{align}
  &\Delta u=0,\ \ (x,y)\in \Omega, \\
  &u(x,y)\mid_{\tau}=f(x,y). \label{su}
\end{align}
Here $\tau$ is a piece of parabola which satisfies $\tau\subset T\subset\subset \Omega$. Let $P$ be a disk with radius $r$, which satisfies $\tau\subset T\subset P\subset \Omega$. By the proof of Lemma \ref{lemmac},
\begin{align}\label{u}
  u(x,y)=\int_{\partial P}\log((x-\zeta_{1})^{2}+(y-\zeta_{2})^{2})\mu(\mathbf{\zeta})ds_{\mathbf{\zeta}},\ \ x\in P.
\end{align}
Discretize the integral of u in (\ref{u}). Apply collocation method, collocate points are $\{(x^{i},y^{i})\}_{i=1}^{I}\in T$, where
$\{(x^{i},y^{i})\}_{m}^{n}\in\tau,\ 1\leq m<n\leq I$, $\{(\zeta_{1}^{j},\zeta_{2}^{j})\}_{j=1}^{J}\in \partial P$. Notice that different choices of $P$ and $(\zeta_{1}^{j},\zeta_{2}^{j})$ will result in different accuracy of the reconstruction. More details, which will not be repeated here, are discussed in \cite{bib:8}.
By (\ref{su}) and (\ref{u})
\begin{align}\label{du}
  \begin{split}
    f(x^{i},y^{i})&=u(x^{i},y^{i})=\int_{\partial P}\log((x^{i}-\zeta_{1})^{2}+(y^{i}-\zeta_{2})^{2})\mu(\mathbf{\zeta})ds_{\mathbf{\zeta}} \\
    &\approx \sum\limits_{j=1}^{J}\mu_{j}\log((x^{i}-\zeta_{1}^{j})^{2}+(y^{i}-\zeta_{2}^{j})^{2})\times \frac{2\pi r}{J},\ \ i=m,m+1,\ldots,n.
  \end{split}
\end{align}
Thus we have
\begin{align}\label{fdu}
  \mathbf{f}=\mathbf{K}\mathbf{\mu},
\end{align}
where $\mathbf{f}=\{f(x^{i},y^{i})\}_{i=1}^{I}$, $\mathbf{K}=(K_{ij})$ is an $I\times J-$order matrix, where $K_{ij}=\log((x^{i}-\zeta_{1}^{j})^{2}+(y^{i}-\zeta_{2}^{j})^{2})\times \frac{2\pi r}{J}$, and $\mathbf{\mu}$ is a $J\times 1$ vector.Then the value of $\mathbf{\mu}=\{\mu_{j}\}_{j=1}^{J}$ can be calculated from (\ref{fdu}).

Note that the collocation linear equation (\ref{du}) is an ill-posed problem, which brings instability due to observation errors. It is generally known that regularization is a method that can improve the stability of ill-posed problems. A deterministic regularization based on Tikhonov regularization in \cite{bib:14} can be adopted here to reduce the instability. Since the conditional stability estimate of $u$ has been obtained in Section \ref{main}, reasonable constraints can be obtained by using the deterministic regularization and also a basis for the selection of regularization parameters can be provided. Then for
\begin{align}
 F=\|K\mathbf{\mu}-\mathbf{f}\|^{2}+\alpha\|\mathbf{\mu}\|^{2},
\end{align}
an optimal $\mu$ can minimize $F$, that is
\begin{align}
  \min\limits_{\mathbf{\mu}} F=\|K\mathbf{\mu}-\mathbf{f}\|^{2}+\alpha\|\mathbf{\mu}\|^{2},
\end{align}
Then
\begin{align}
  \mathbf{\mu}=(\alpha\mathbf{I}+K^{*}K)^{-1}K^{*}\mathbf{f}.
\end{align}
Here $\alpha$ is a regularization parameter, which can be chosen as a number with the same order with  observation error \cite{bib:14}.

After $\{\mu_{j}\}_{j=1}^{J}$ being acquired, $u\mid_{T}$ can be approximated by
\begin{align}
  u(x,y)=\sum\limits_{j=1}^{J}\mu_{j}\log((x-\zeta_{1}^{j})^{2}+(y-\zeta_{2}^{j})^{2})\frac{2\pi r}{J},\ \ (x,y)\in T,
\end{align}
the value of $u$ on $T$ can be obtained like a cork.

Here are two applications used to illustrate the conditional stability estimation. The settings of the numerical cases are summarized in Table \ref{case-list}.

\subsection{Parabolic curve}

Assume that $T=\{(x,y)\mid y=-0.5+2x^{2},x\in(-0.5,0.6)\}$, and $f(x,y)=e^{-2x}\cos(2y)$, then the collocation points on $T$ are $\{(x^{i},y^{i})=(-0.5+\frac{0.35\pi i}{I},-0.5+2(-0.5+\frac{0.35\pi i}{I})^{2})\}_{i=1}^{I},I=180$. $\mathbf{x_{0}}$ is selected as $\mathbf{x_{0}}=(0,0)$ in calculations. Construct $P$ as a disk with center $\mathbf{x_{0}}$ and radius $r=1.1$. Then the collocation points on $\partial P$ are $\{(\zeta_{1}^{j},\zeta_{2}^{j})=(r\cos\frac{2\pi j}{J},r\sin\frac{2\pi j}{J})\}_{j=1}^{J},\ J=40$.

Consider the cases of one or two fixed segments, i.e. $\tau$ is taken as a segment or the union of two segments on $T$. Then consider the $\tau$ in different domain of $x$.
\begin{table}
\begin{center}
\begin{tabular}{|c|c|c|}
  \hline
  case & domain of $x$ & Number of collocation points \\
  \hline
  a & $(-0.5,-0.3)$              & $1-30$ \\
  b & $(-0.5,-0.1)$              & $1-60$ \\
  c & $(-0.2,0)$                 & $51-80$ \\
  d & $(-0.5,-0.4)\cup(0.5,0.6)$ & $1-15,166-180$ \\
  e & $(-0.3,-0.2)\cup(0.3,0.4)$ & $31-45,136-150$ \\
  f & $(-0.2,-0.1)\cup(0.0,0.1)$ & $46-60,81-95$ \\
  \hline
\end{tabular}
\end{center}
    \caption{List of numerical cases.}
    \label{case-list}
\end{table}

\begin{figure}[!ht]
\centering
\vspace{-0.35cm}
  \subfigure[$1-30$]{
    \label{fig:subfig:a1} 
    \includegraphics[width=2.6in]{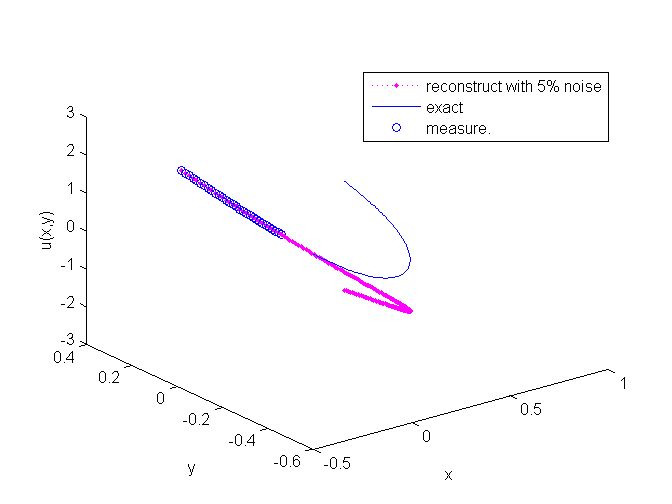}
  }
  \subfigure[$1-60$]{
    \label{fig:subfig:b1} 
    \includegraphics[width=2.6in]{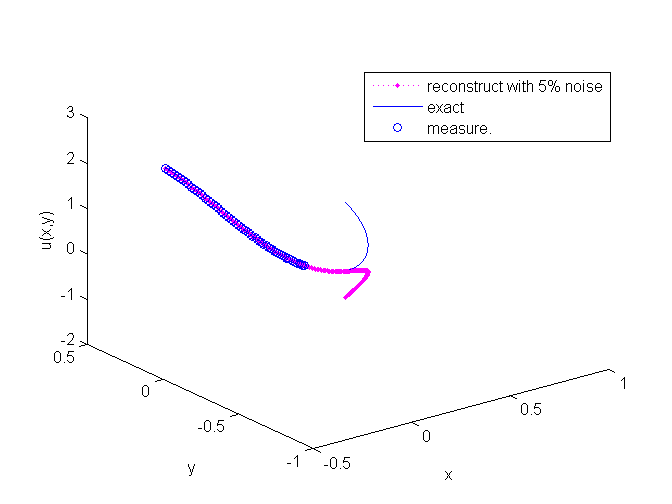}
  }
  \subfigure[$51-80$]{
    \label{fig:subfig:c1} 
    \includegraphics[width=2.6in]{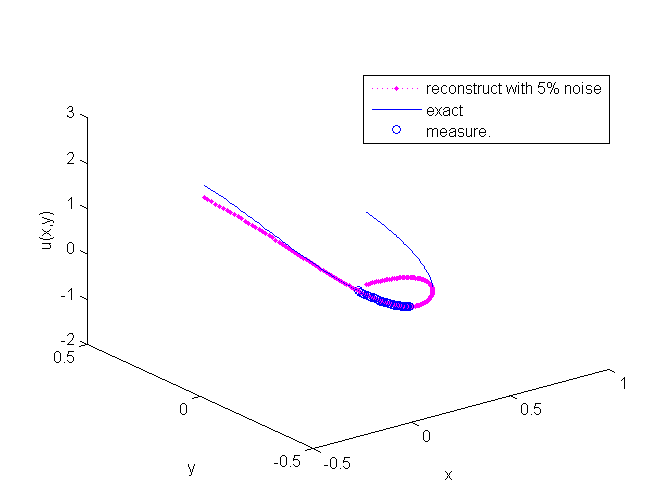}
  }
  \subfigure[$1-15,166-180$]{
    \label{fig:subfig:d1} 
    \includegraphics[width=2.6in]{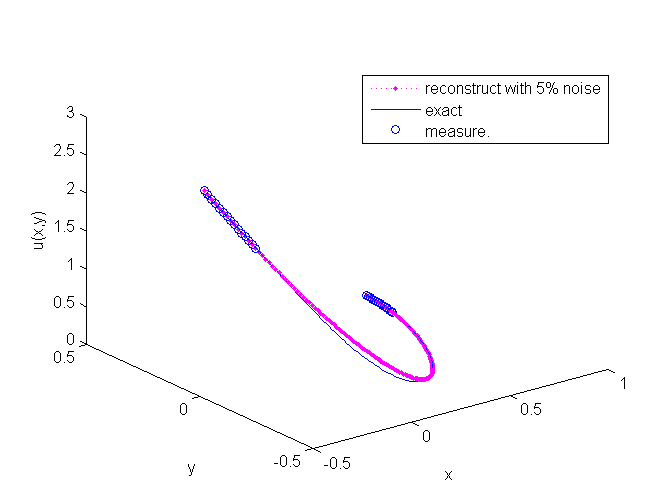}
  }
  \subfigure[$31-45,136-150$]{
    \label{fig:subfig:f1} 
    \includegraphics[width=2.6in]{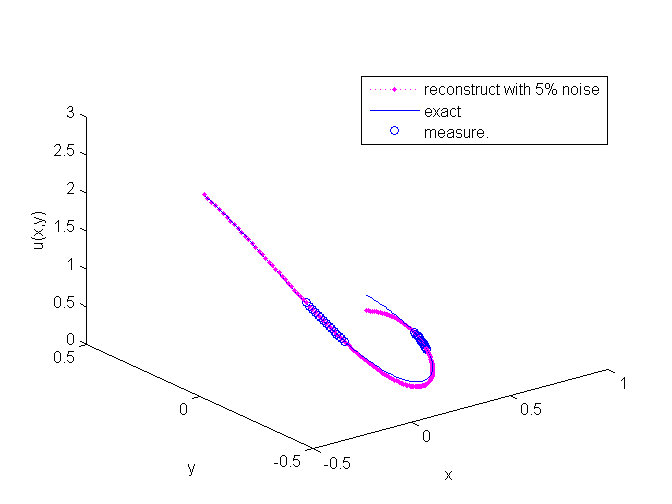}
  }
  \subfigure[$46-60,81-95$]{
    \label{fig:subfig:g1} 
    \includegraphics[width=2.6in]{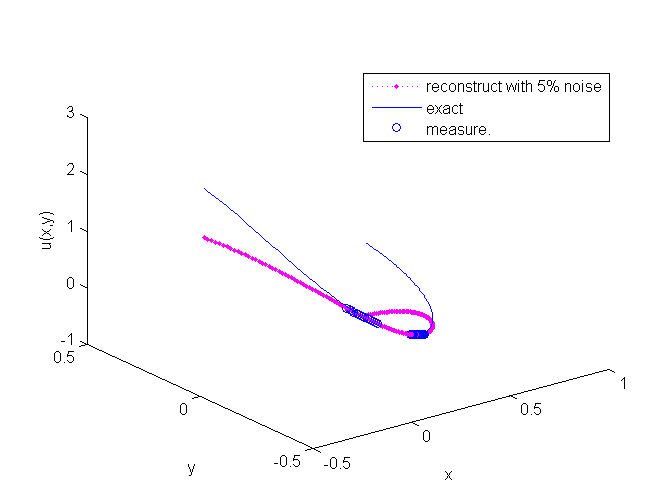}
  }
  \vspace{-0.4cm}
  \caption{Unique continuation on parabolic equation for harmonic function: the capital of each subfigure means the collocation points on $\tau$}
  \label{parabola}
\end{figure}

A point-by-point observation error of $1\%\sim5\%$ is added to $f(x,y)|_{\tau}$ in the above cases respectively. The reconstructed $u(x,y)|_{T}$ is shown in Fig. \ref{parabola}. It can be seen that the value of $u(x,y)$ on $T\backslash\tau$ can be obtained from the value of $u(x,y)$ on $\tau$, which illustrates that it is uniformly continuable.

\begin{figure}[!ht]
\centering
\vspace{-0.35cm}
  \subfigure[$1-30$]{
    \label{fig:subfig:a2} 
    \includegraphics[width=2.6in]{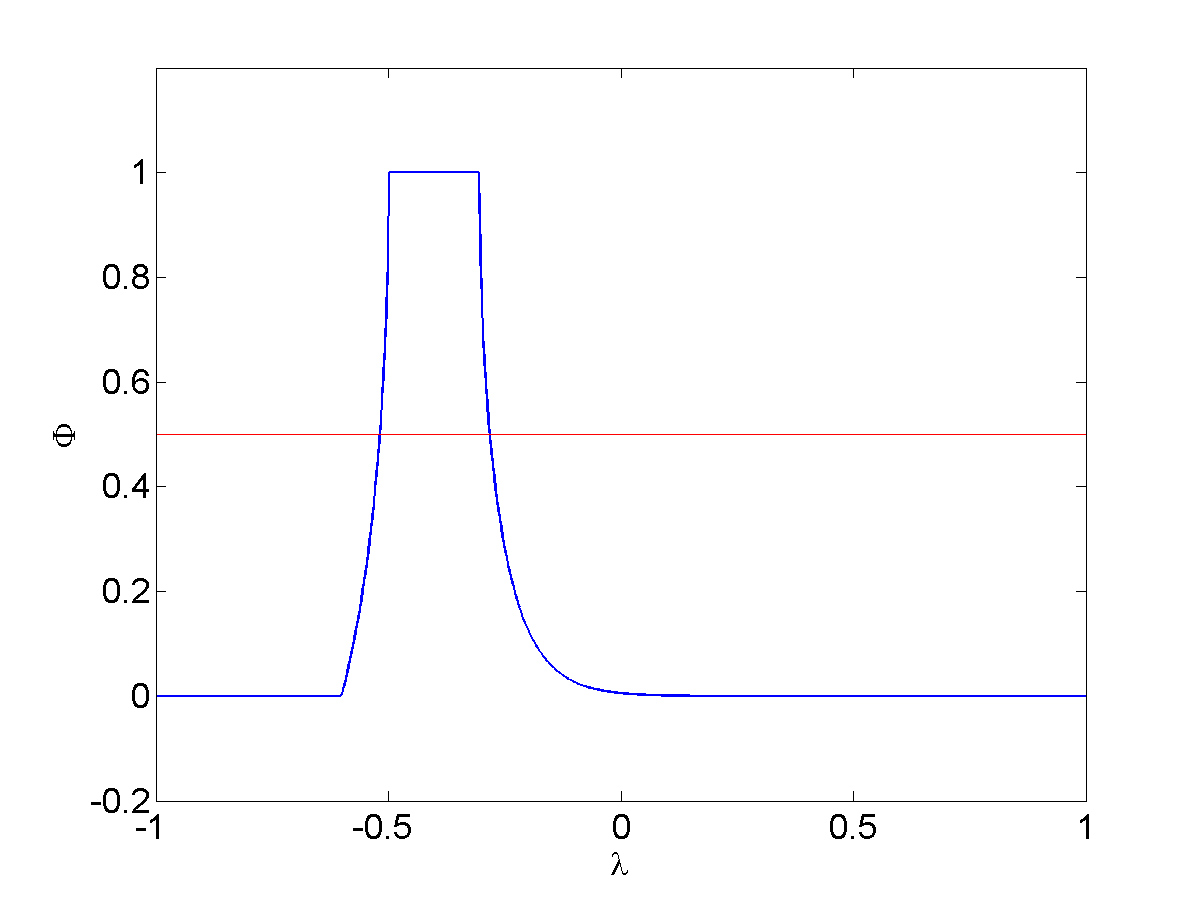}
  }
  \subfigure[$1-60$]{
    \label{fig:subfig:b2} 
    \includegraphics[width=2.6in]{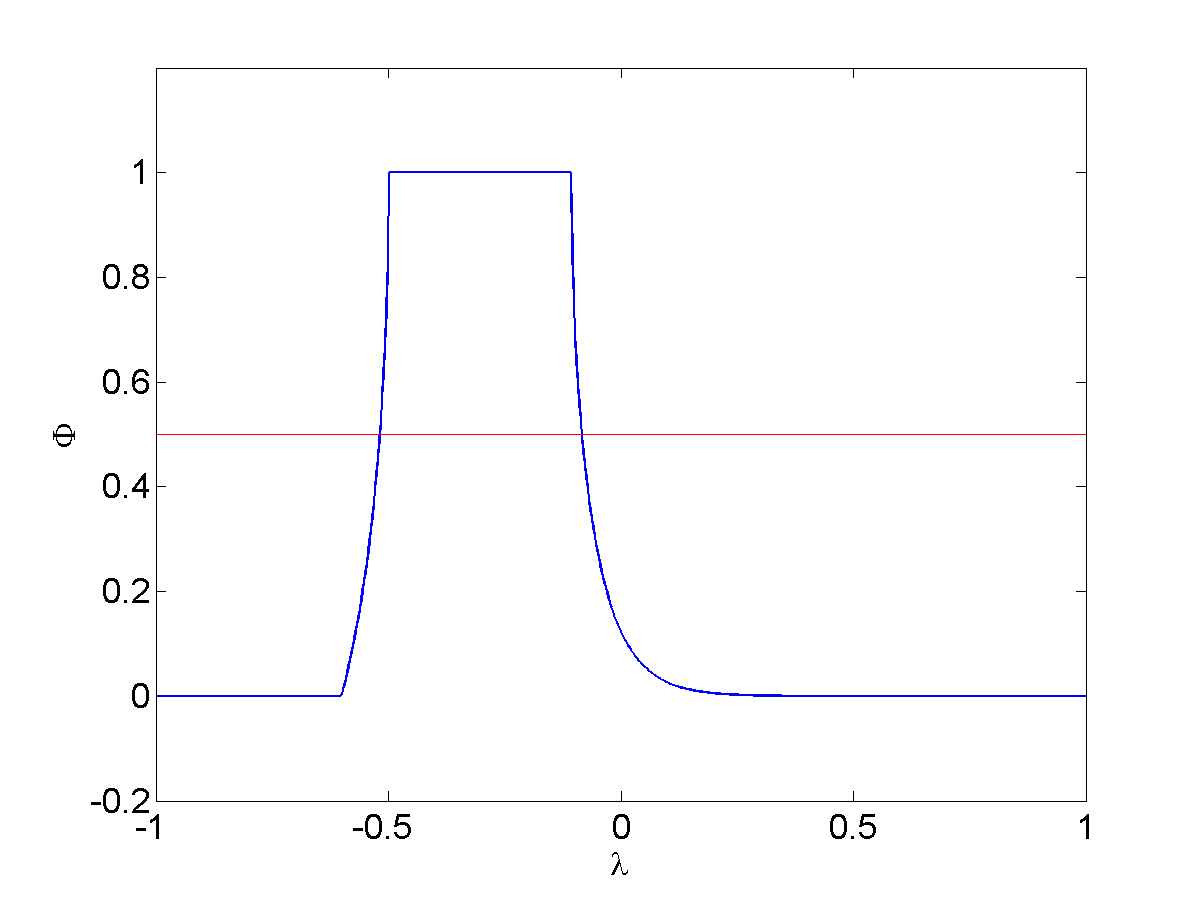}
  }
  \subfigure[$51-80$]{
    \label{fig:subfig:c2} 
    \includegraphics[width=2.6in]{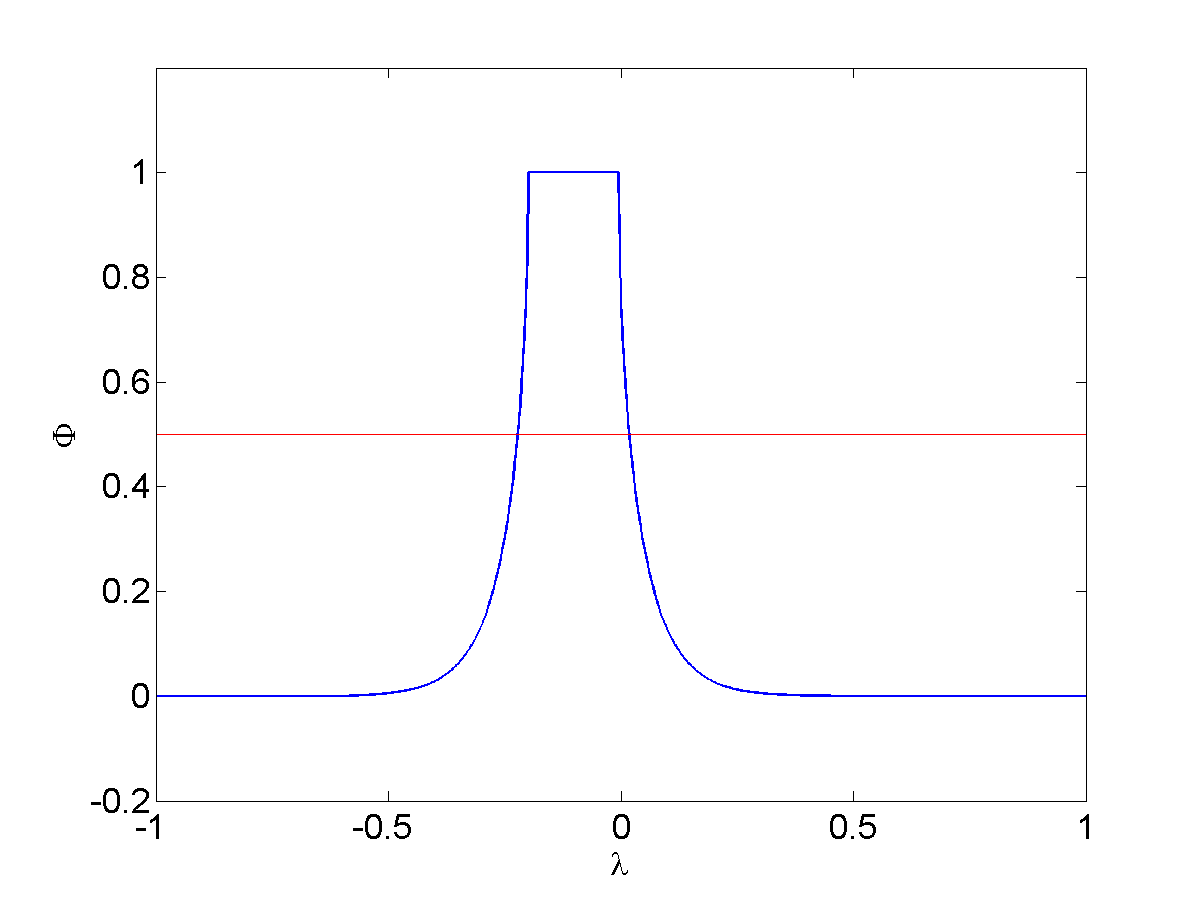}
  }
  \subfigure[$1-15,166-180$]{
    \label{fig:subfig:d2} 
    \includegraphics[width=2.6in]{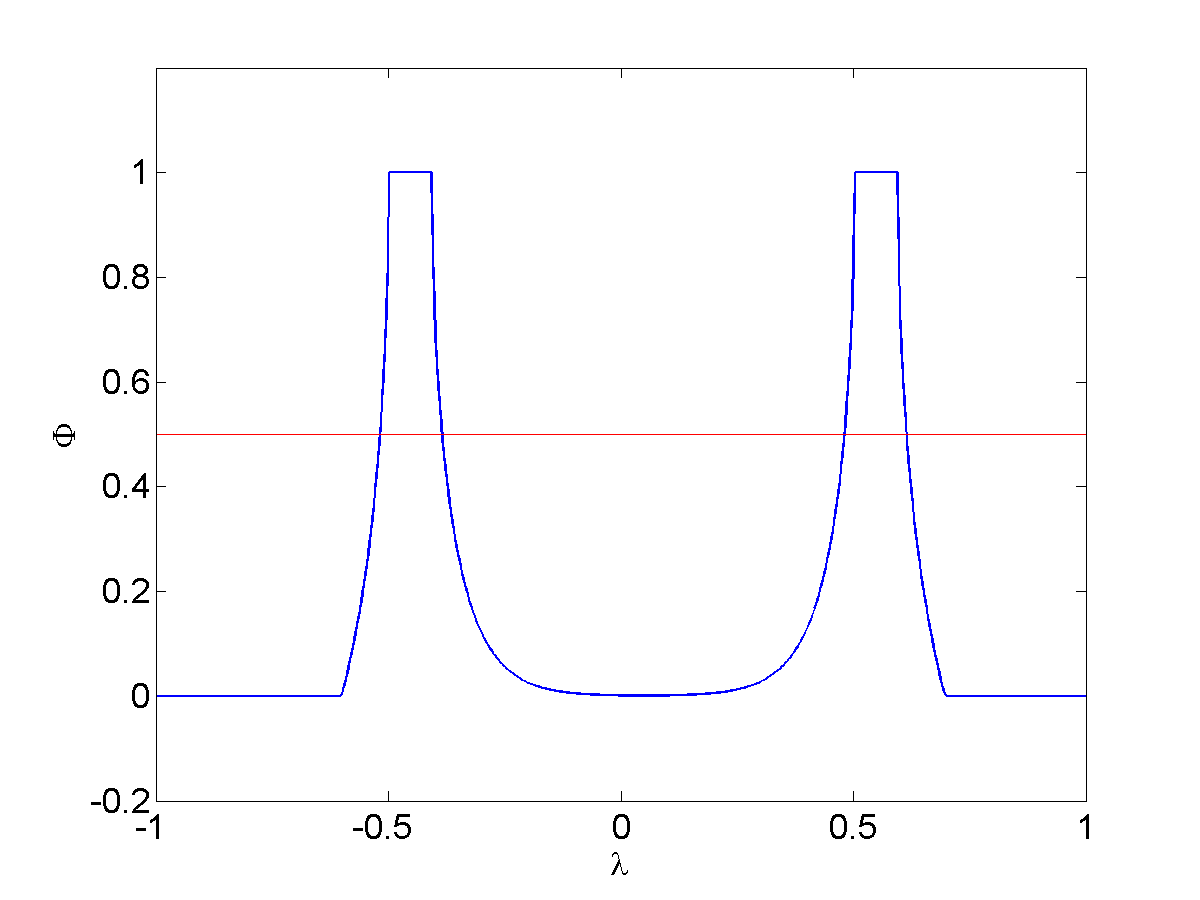}
  }
  \subfigure[$31-45,136-150$]{
    \label{fig:subfig:f2} 
    \includegraphics[width=2.6in]{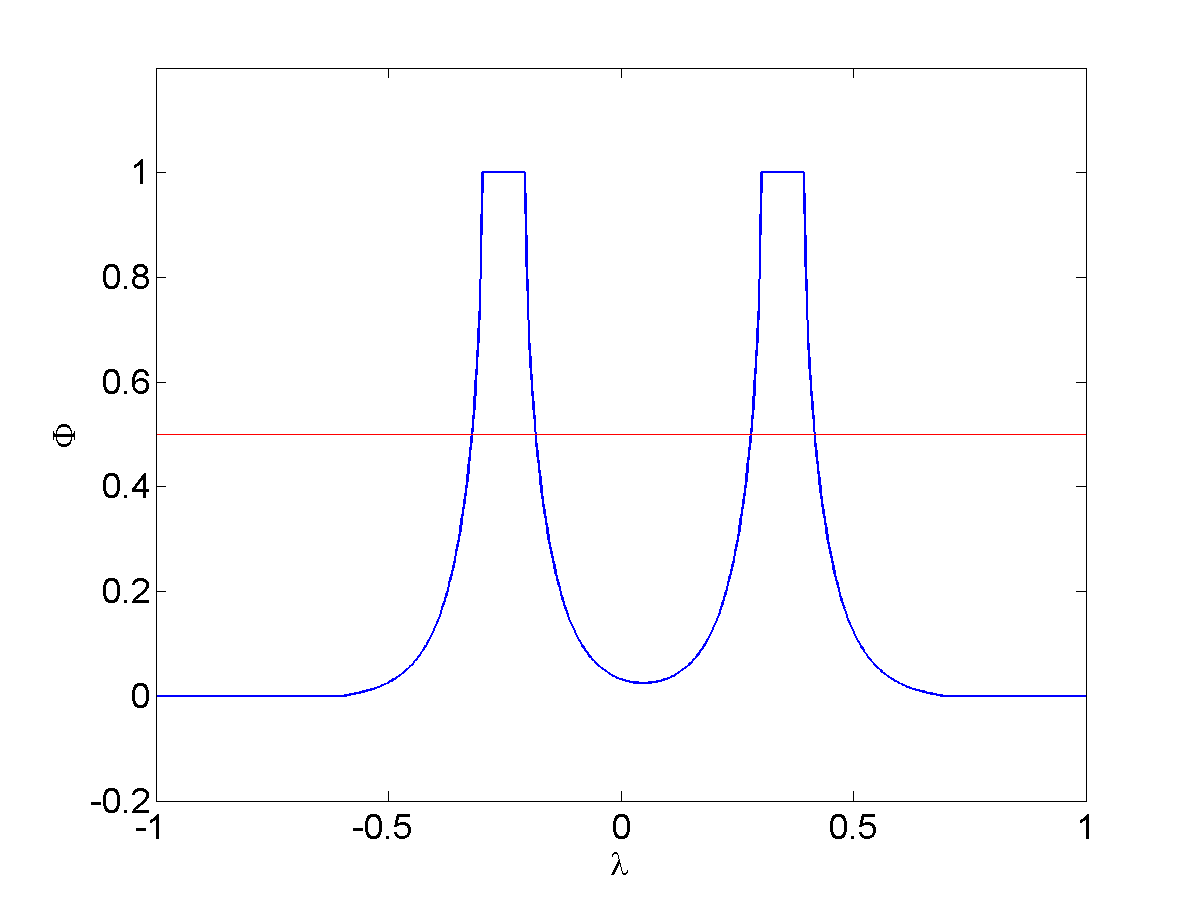}
  }
  \subfigure[$46-60,81-95$]{
    \label{fig:subfig:g2} 
    \includegraphics[width=2.6in]{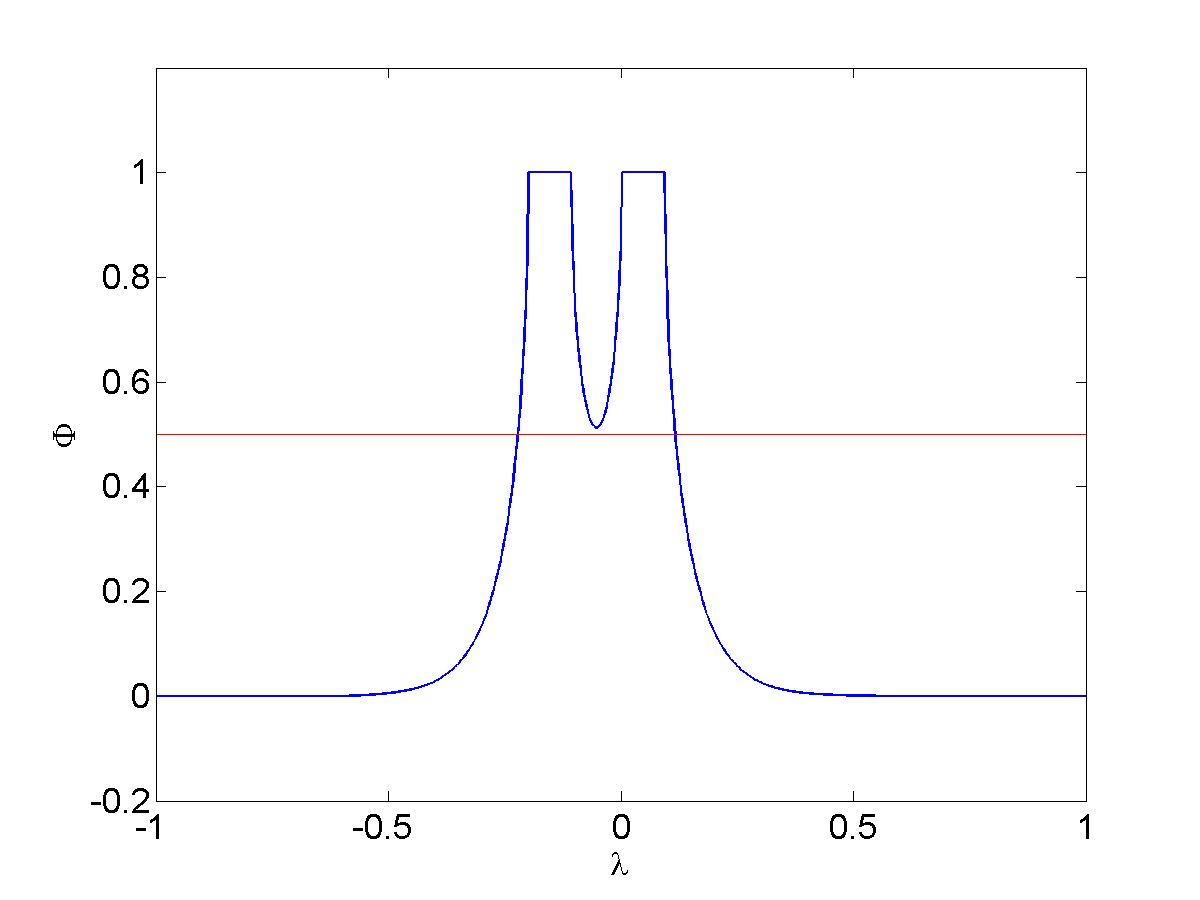}
  }
  \vspace{-0.4cm}
  \caption{Harmonic measure: the capital of each subfigure means the collocation points on $\tau$}
  \label{parabolahm}
\end{figure}
From Fig. \ref{parabola}, comparing Fig. \ref{fig:subfig:a1} and Fig. \ref{fig:subfig:b1}, more fixed points have no significant refluence on control of the value of $u$. When the one end is fixed, the reconstruction of $u$ turns to be dissatisfactory when it extents over the vertex of $T$. If the points near to vertex of $T$ are fixed, the reconstruction of $u$ in Fig. \ref{fig:subfig:c1} is better than that in Fig. \ref{fig:subfig:a1} in total, however, the value is still inaccurate on the both end of $T$. Thus, dividing the fixed points in two sections is considered. It can be seen in the rest figures in Fig. \ref{parabola} that when the two sections are far enough, the whole line can be better controlled, but the line between two sections may not good enough. The reason can be seen from harmonic measure.

Remark that the harmonic measure $\varphi$ decides the degree of control for u on T, from Fig. \ref{parabolahm}, if we define that the points with $\varphi>\frac{1}{2}$ are good-controlled points, then we can see that, the good-controlled points only appear when nearing to the controlled points. Consider the situation of two fixed sections, the line of $\varphi$ between two fixed lines can be controlled better. However, if two sections are far away enough, the line in the middle between them still cannot be controlled well.

\begin{table}[!ht]
    \centering
    \subtable[observation error of $1\%$]{
        \begin{tabular}{| c | c | c | c |}
            \hline
                  domain of controlled $x$ & number of controlled points & $\|f-u\|_{\tau}$ & $\|f-u\|_{T}$ \\
            \hline
                  $(-0.5,0.3)$ & 30 &  0.0116 & 4.8247  \\
                  $(-0.2,0)$ & 30 & 0.0038 & 2.0879 \\
                  $(-0.5,-0.4)\cup(0.5,0.6)$ & 30 & 0.0089 & 0.9060 \\
                  $(-0.5,-0.3)\cup(0.4,0.6)$ & 60 & 0.0065 & 0.3274 \\
            \hline
        \end{tabular}
    }
    \subtable[observation error of $5\%$]{
     \begin{tabular}{| c | c | c | c |}
        \hline
              domain of controlled $x$ & number of controlled points & $\|f-u\|_{\tau}$ & $\|f-u\|_{T}$ \\
        \hline
              $(-0.5,-0.3)$ & 30 &  0.0549 & 4.9498  \\
              $(-0.2,0)$ & 30 & 0.0205 & 2.3915 \\
              $(-0.5,-0.4)\cup(0.5,0.6)$ & 30 & 0.0286 & 1.4035 \\
              $(-0.5,-0.3)\cup(0.4,0.6)$ & 60 & 0.0359 & 0.4281 \\
        \hline
    \end{tabular}
    }
    \caption{Error on the parabolic curve}
    \label{ptable}
\end{table}
From Table \ref{ptable}, it can be seen that when the controlled points are same, it is obvious that the smaller observation error, the more accurate the estimation result. When the observation error is controlled at the same level, the accuracy of the estimated results from Case a to Case d is significantly improved. In Case d, when the first 30 points and the last 30 points of the 180 points are controlled, as can be seen from the table, even if there is an error of 5\% in the observed value $u|_{\tau}$, the estimated result $u|_{T}$ is still very close to the accurate value $f$, the error is only of order $10^{-1}$.

\subsection{Hyperbolic curve}

After the discussion on parabolic equation, according to the same method, the unique continuation for harmonic functions on other curves can be easily obtained.

Assume that $T=\{(x,y)\mid \frac{y^{2}}{0.25}-\frac{x^{2}}{0.36}=1,x\in(-0.5,0.6)\}$, then the collocation points on $T$ are $\{(x^{i},y^{i})=(-0.5+\frac{0.35\pi i}{I},\frac{5}{6}\sqrt{(-0.5+\frac{0.35\pi i}{I})^{2}+0.36}\}_{i=1}^{I},I=180$. The values of $f(x,y)$, $\mathbf{x_{0}}$, $P$, $r$ and the collocation points on $\partial P$: $\{(\zeta_{1}^{j},\zeta_{2}^{j})=(r\cos\frac{2\pi j}{J}),r\sin\frac{2\pi j}{J})\}_{j=1}^{J},\ J=40$ are the same selections as before.

Consider the cases of one or two fixed curves again with the same fixed points.

\begin{figure}[!ht]
\vspace{-0.35cm}
  \subfigure[$1-30$]{
    \label{fig:subfig:a3} 
    \includegraphics[width=2.6in]{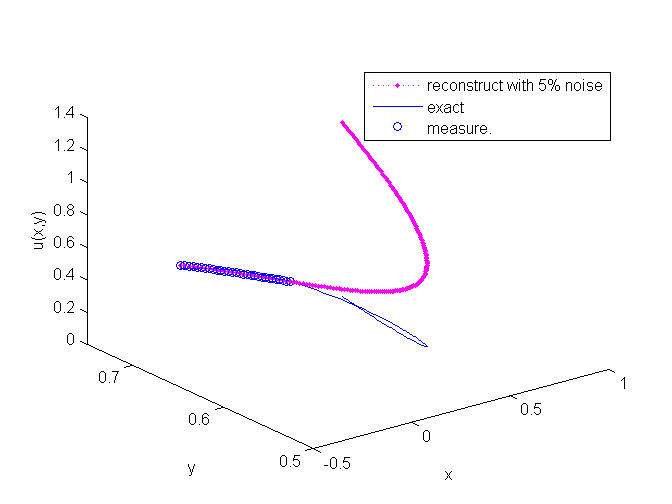}
  }
  \subfigure[$1-60$]{
    \label{fig:subfig:b3} 
    \includegraphics[width=2.6in]{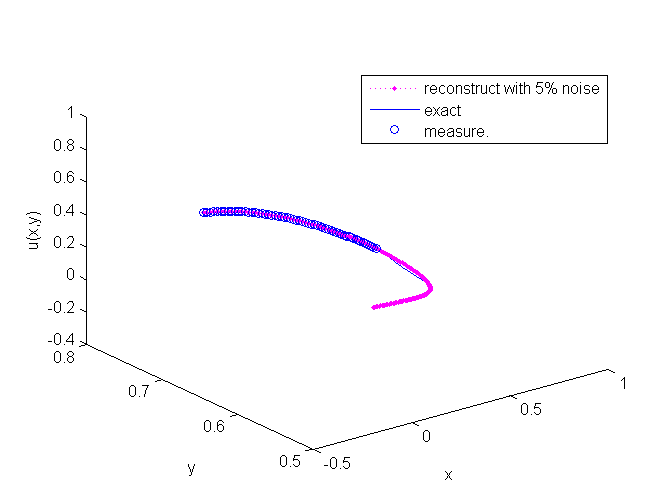}
  }
  \subfigure[$51-80$]{
    \label{fig:subfig:c3} 
    \includegraphics[width=2.6in]{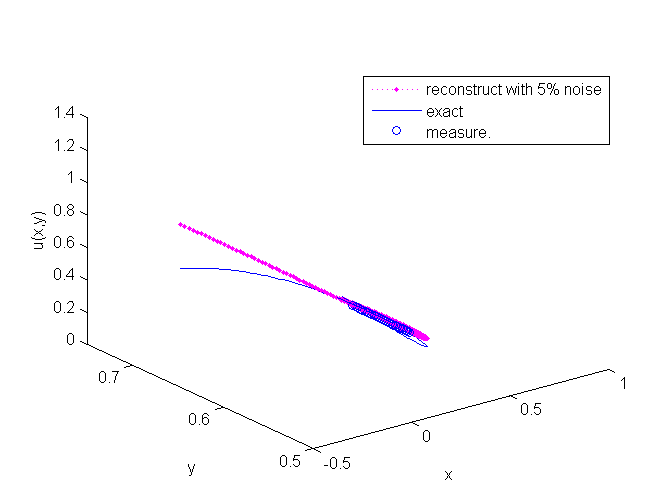}
  }
  \subfigure[$1-15,166-180$]{
    \label{fig:subfig:d3} 
    \includegraphics[width=2.6in]{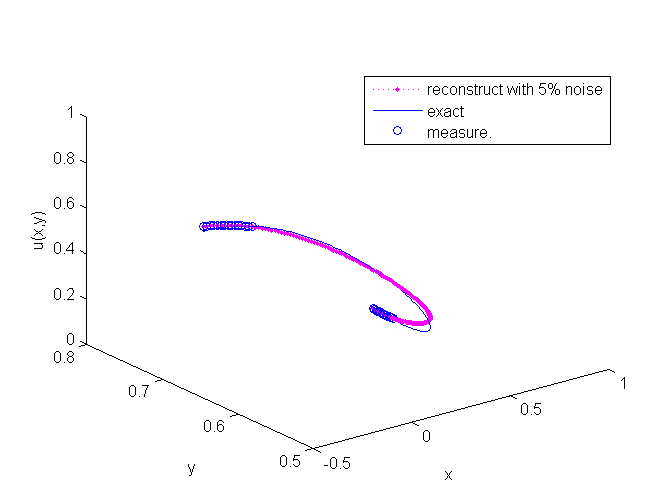}
  }
  \subfigure[$31-45,136-150$]{
    \label{fig:subfig:f3} 
    \includegraphics[width=2.6in]{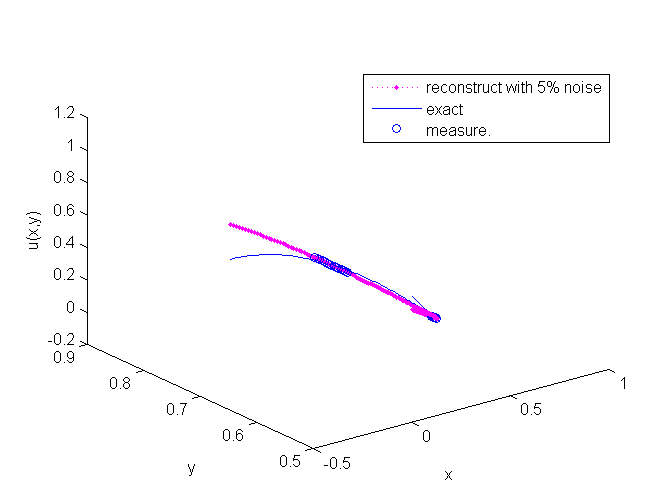}
  }
  \subfigure[$46-60,81-95$]{
    \label{fig:subfig:g3} 
    \includegraphics[width=2.6in]{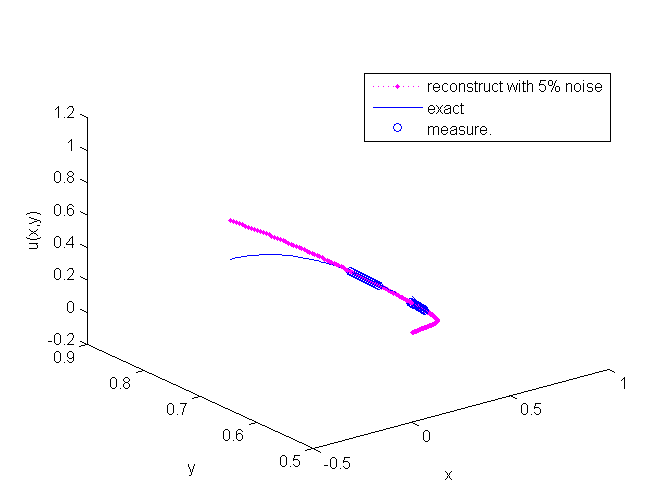}
  }
  \vspace{-0.4cm}
  \caption{Unique continuation on hyperbolic equation for harmonic function: the capital of each subfigure means the collocation points on $\tau$}
  \label{logarithm} 
\end{figure}

Similarly, a point-by-point observation error of $1\%\sim5\%$ is added to $f(x,y)|_{\tau}$ in the above four cases respectively. The reconstructed $u(x,y)|_{T}$ is shown in Fig. \ref{logarithm}. It can be seen that the value of $u(x,y)$ on $T\backslash\tau$ can be obtained from the value of $u(x,y)$ on $\tau$, which can be proved uniformly continuable. And similar results hold on hyperbolic curve with parabola curve from Fig. \ref{logarithm}.

\begin{table}
    \centering
    \subtable[observation error of $1\%$]{
        \begin{tabular}{| c | c | c | c |}
            \hline
                  domain of controlled $x$ & number of controlled points & $\|f-u\|_{\tau}$ & $\|f-u\|_{T}$ \\
            \hline
                  $(-0.5,0.3)$ & 30 &   0.0031 & 1.3896  \\
                  $(-0.2,0)$ & 30 & 0.0031 & 0.3103 \\
                  $(-0.5,-0.4)\cup(0.5,0.6)$ & 30 & 0.0023 & 0.0850 \\
                  $(-0.5,-0.3)\cup(0.4,0.6)$ & 60 & 0.0021 &   0.0572 \\
            \hline
        \end{tabular}
    }
    \subtable[observation error of $5\%$]{
     \begin{tabular}{| c | c | c | c |}
        \hline
              domain of controlled $x$ & number of controlled points & $\|f-u\|_{\tau}$ & $\|f-u\|_{T}$ \\
        \hline
              $(-0.5,0.3)$ & 30 &  0.0233 & 1.6567  \\
              $(-0.2,0)$ & 30 & 0.0159 & 0.3509 \\
              $(-0.5,-0.4)\cup(0.5,0.6)$ & 30 & 0.0125 & 0.1668 \\
              $(-0.5,-0.3)\cup(0.4,0.6)$ & 60 &  0.0111 & 0.0763 \\
        \hline
    \end{tabular}
    }
    \caption{Error on the hyperbolic curve}
    \label{ltable}
\end{table}
Table \ref{ltable} shows the similar conclusion as the conclusion from Table \ref{ptable}, which illustrates that the conditional stability estimate in (\ref{ux}) is valid.

The specific results of multiple sections will be discussed in subsequent papers.

\section{Conclusion}\label{conclusion}

This paper obtains the unique continuation on quadratic curves for harmonic functions. Similar to the results on a straight line for harmonic functions, the unique continuation and the conditional stability are proved. A difference is that the complex extension on the quadratic curves is more complicated. Due to the complexity of the complex extension of the function on the quadratic curve, it is necessary to consider the selection of the complex extension region and the harmonic measure.

In this paper, the unique continuation of two types of curves (parabola and hyperbola) are calculated numerically by means of collocation method and Tikhonov regularization. The calculation takes into account the length, position and number of segments of the controlled curves. Through comparison, it is found that the lengths of control range $\tau$ has no significant results on estimation of $u\mid_{T}$, but the multi-stage control can achieve more accurate results. These numerical applications are consistent with the theoretical results in this article and \cite{bib:3}. Thus in applications, when data on a quadratic curve need to be measured, according to the description of this paper, only a part of data, instead of all the data  needs to be measured directly, which reduces the measurement cost significantly.

\bigskip

{\bf Acknowledgement}\ \ This work was supported by National Science Foundation Committee(No. 11971121).

\end{document}